\documentclass[runningheads]{llncs}

\usepackage{graphicx}
\usepackage{amssymb,mathrsfs}
\usepackage{amsmath}
\usepackage{stmaryrd}
\usepackage{tikz}
\usepackage{color}
\usepackage{enumitem}
\usepackage{tablefootnote}
\usepackage[titletoc,title]{appendix}
 \usepackage[title]{appendix}
\usepackage[all,pdf]{xy}
\usepackage{bussproofs}

\usepackage{multirow}

\newcommand{\D}{\Diamond}
\newcommand{\B}{\Box}
\newcommand{\eq}{\leftrightarrow}

\renewcommand{\phi}{\varphi}

\newcommand{\mf}{\mathfrak}

\newcommand{\mc}{\mathcal}

\newcommand{\sm}{\setminus}
\newcommand{\sub}{\subseteq}

\newcommand{\cset}[1]{{\{ #1 \}}}
\newcommand{\tup}[1]{{\langle #1 \rangle}}
\newcommand{\bd}{\blacklozenge}
\newcommand{\bb}{\blacksquare}

\renewcommand{\D}{\Diamond}
\newcommand{\M}{\mf M}

\newcommand{\LHS}{\mathsf{LHS}}
\renewcommand{\L}{\mc{L}}

\newcommand{\md}{\models}

\newcommand{\N}{\mathfrak{N}}

\newcommand{\bis}{\mathrel{\mathchoice%
{\raisebox{.3ex}{$\,
  \underline{\makebox[.7em]{$\leftrightarrow$}}\,$}}%
{\raisebox{.3ex}{$\,
  \underline{\makebox[.7em]{$\leftrightarrow$}}\,$}}%
{\raisebox{.2ex}{$\,
  \underline{\makebox[.5em]{\scriptsize$\leftrightarrow$}}\,$}}%
{\raisebox{.2ex}{$\,
  \underline{\makebox[.5em]{\scriptsize$\leftrightarrow$}}\,$}}}}

\begin{document}
\title{Logic of the Hide and Seek Game: Characterization, Axiomatization, Decidability}
%
\titlerunning{LHS: Characterization, Axiomatization, Decidability}
%

\author{~}
\institute{~}

\author{Qian Chen\inst{1} \and Dazhu Li\inst{2,3}}
 
\authorrunning{Chen $\&$ Li}
 
\institute{Department of Philosophy,
Tsinghua University, Beijing, China\\
\email{chenq21@mails.tsinghua.edu.cn}  \and
Institute of Philosophy, Chinese Academy of Sciences, Beijing, China 
\and
Department of Philosophy, University of Chinese Academy of Sciences\\
\email{lidazhu@ucas.ac.cn}}

\maketitle               


\begin{abstract}
The logic of the hide and seek game $\LHS$ was proposed to reason about search missions and interactions between agents in pursuit-evasion environments. 
As proved in \cite{Li2021,LHS-journal}, having an equality constant in the language of $\LHS$ drastically increases its  computational complexity: the satisfiability problem for $\LHS$ with multiple relations is undecidable. In this work, we improve the existing result by showing that $\LHS$ with a single relation is undecidable. With the findings of \cite{Li2021,LHS-journal}, we provide a van Benthem style characterization theorem for the expressive power of the logic. Finally, by `splitting' the language of $\LHS^-$, a crucial fragment of $\LHS$ without the equality constant, into two `isolated parts', we provide a complete Hilbert style proof system for $\LHS^-$ and prove that its satisfiability problem is decidable, whose proofs would indicate significant differences between the proposals of $\LHS^-$ and of ordinary product logics. Although $\LHS$ and $\LHS^-$ are frameworks for interactions of 2 agents, all results in the article can be easily transferred to their generalizations for settings with any $n>2$ agents. 
\end{abstract}

\vspace{-.7cm}


\section{Introduction}\label{sec:introdunction}
The logic $\LHS$ of the hide and seek game was introduced in \cite{graph-game-logic-design} that promotes a study of graph game design in tandem with matching modal logics, and then was probed in \cite{Li2021} and its further extension \cite{LHS-journal}. The logic provides us with a platform to reason about search problems and  interactions between agents with entangled goals, as in the case of the hide and seek game \cite{graph-game-logic-design} (or the game of cops and robber \cite{cop-robber}): in a fixed graph, two players Hider and Seeker take turns to move to a successor of their own positions, and Seeker tries to move to the same position with Hider while Hider aims to avoid Seeker. 

To describe the game, the language of $\LHS$ contains two modalities for the movements of the two players and a constant $I$ expressing that the positions of Hider and Seeker are the same. Semantically, models for $\LHS$ are the same as relational models for basic modal logic \cite{BRV2001}, while formulas are evaluated at {\em two} states, which intuitively represent the positions of the two players.

In addition to the applications to the graph games, $\LHS$ is also of interest from  other perspectives. One of them is that the framework links up the study of graph game logics with many other important fields: as illustrated in \cite{Li2021,LHS-journal}, $\LHS$ and {\em its fragment $\LHS^-$ without the constant $I$} have close connection with {\em product logics}, including $\mathsf{K\times K}$ \cite{many-dimensional-book} and its extension $\mathsf{K\times^{\delta} K}$ with {\em a diagonal constant $\delta$} \cite{product-undecidable,Kikot2010,diagonal-fmp}; the framework $\LHS$ is highly relevant to {\em cylindric modal logics} that also contain constants for equality \cite{Yde-thesis}; and both $\mathsf{K\times^{\delta} K}$ and cylindric modal logics in turn provide a link between $\LHS$ and {\em cylindric algebra} proposed in \cite{cylindric-algebra}. Moreover, the framework $\LHS$  provides an instance showing how an innocent looking proposal $I$ for equality can drastically improve the computational complexity of the logic: as proved in \cite{Li2021,LHS-journal}, the satisfiability problem for $\LHS$ with multiple binary relations is undecidable.

In this work, we will explore the further properties of $\LHS$ and $\LHS^-$. First, we  improve the existing undecidability result for $\LHS$ with multiple relations and show that $\LHS$ with a single relation is undecidable (Section \ref{sec:undecidability}). Then, based on the notions of first-order translation and bisimulations for $\LHS$ given in \cite{LHS-journal}, we develop a van Benthem style characterization theorem for the expressiveness of $\LHS$ (Section \ref{sec:expressiveness}). Next, for  $\LHS^-$, we develop a complete Hilbert style calculus and show that its satisfiability problem is decidable (Section \ref{sec:axiomatization}), and our proofs would indicate important differences between the proposals of $\LHS^-$ and $\mathsf{K\times K}$. Also, we discuss related work and point out a few lines for further study (Section \ref{sec:conclusion}). It is instructive to notice that although $\LHS$ and $\LHS^-$ are frameworks for the hide and seek game with 2 players, all these results can be transferred to the logics generalizing $\LHS$ and $\LHS^-$ for the settings with $n>2$ players, but we stick to discussing the systems $\LHS$ and $\LHS^-$  for simplicity.

\vspace{3mm}


\section{Basics of the Logic of the Hide and Seek Game}\label{sec:logic-design}

We start by concisely introducing the basics of $\mathsf{LHS}$, including its language and semantics, and providing preliminary observations on its properties. 

    \begin{definition}
        Let $\mathsf{L}=\cset{p^l_i:i\in\mathbb{N}}$ and $\mathsf{R}=\cset{p^r_i:i\in\mathbb{N}}$ be two disjoint countable sets of propositional variables. The {\em language} $\L$ of $\mathsf{LHS}$ is given by:
        \[
            \L\ni \phi ::= p^l \mid p^r \mid I \mid \neg\phi \mid \phi\wedge\phi \mid \B\phi \mid \bb\phi,
        \]
        where $p^l\in \mathsf{L}$ and $p^r\in \mathsf{R}$.  
    \end{definition}
    
Abbreviations $\top,\bot,\vee,\to$ are as usual, and we use $\D,\bd$ for the dual operators of $\B$ and $\bb$ respectively. For convenience, we call $\B$ and $\D$ `{\em white modalities}' and call $\bb$ and $\bd$ `{\em black modalities}'. Also, the notion of {\em subformulas} is as usual, and for any $\varphi\in\L$, we employ $\mathsf{Sub}(\phi)$ for {\em the set of subformulas of  $\phi$}. In what follows, we use $\mathcal{L}^-$ for the part of $\mathcal{L}$ without $I$, which is the language for $\mathsf{LHS}^{-}$.

A {\em frame} is a tuple  $F=(W,R)$ such that $W$ is a non-empty set of states and $R\sub W\times W$ is a binary relation on $W$. A {\em model} $\M=(W,R,V)$ equips a frame with a valuation function $V:\mathsf{L}\cup\mathsf{R}\to\mathcal{P}(W)$.\footnote{For any set $A$, we use $\mathcal{P}(A)$ for its {\em power set}.} For any  $s,t\in W$, we call $\tup{\M,s,t}$ a {\em pointed $\LHS$-model}. For simplicity, we usually write $\M,s,t$ for it. 
For each $w\in W$ and $U\sub W$, we define $R(w)=\cset{v\in W: Rwv}$ and $R(U)=\bigcup_{u\in U}R(u)$.

    \begin{definition}
Let $\M = (W, R,V)$ be a model and $s,t\in W$. {\em Truth of formulas} $\varphi\in\mathcal{L}$ at $\tup{\M,s,t}$, written as $\M,s,t\md\varphi$, is defined recursively as follows:
  \begin{center}
\begin{tabular}{r@{\quad $\Leftrightarrow$\quad}l}
$\M,s,t\md p^l$ & $s\in V(p^l)$  \\
$\M,s,t\md p^r$ & $t\in V(p^r)$  \\
$\M,s,t\md I$& $s=t$\\
$\M,s,t\md \neg\varphi$ & $\M,s,t\not\md\varphi$  \\
$\M,s,t\md \varphi\wedge\psi$ & $\M,s,t\md\varphi$ and $\M,s,t\md\psi$  \\
$\M,s,t\md\B\varphi$& $\M,s',t\md \varphi$ for all $s'\in R(s)$ \\
$\M,s,t\md\bb\varphi$& $\M,s,t'\md\varphi$ for all $t'\in R(t)$ \\
\end{tabular}
\end{center}
\end{definition}

Notions of {\em satisfiability}, {\em validity} and {\em logical consequence} are defined in the usual manner. Let $\mathsf{LHS}$ denote the set of all valid formulas. Let $\M=(W,R,V)$ be a model and $U\sub W$. We say {\em a model $\M'=(W',R',V')$ is generated from $\M=(W,R,V)$ by $U$}, if $\M'$ is the smallest model satisfying the following: $U\sub W'$, $R(W')\sub W'$, $R'=R\cap (W'\times W')$, and for each $p\in \mathsf{L}\cup \mathsf{R}$, $V'(p)= V(p)\cap W'$.

\begin{proposition}\label{prop:generated-model}
    Let $\M'=(W',R',V')$ be a model generated from $\M=(W,R,V)$ by $\{s,t\}\subseteq W$. For any formula $\varphi\in\mathcal{L}$, it holds that: 
        $\M,s,t\md \varphi$ iff $\M',s,t\md \varphi$.
\end{proposition}

\begin{proof}
    It goes by induction on formulas. We omit the details to save space.\qed
\end{proof}

\section{Undecidability of $\mathsf{LHS}$}\label{sec:undecidability}

As stated, \cite{Li2021,LHS-journal} proved that the satisfiability problem for $\mathsf{LHS}$ with multiple binary relations is undecidable. In this part, we show that $\mathsf{LHS}$ with a single relation is undecidable as well, which is an improvement of the existing result.

\begin{theorem}
The satisfiability problem for $\mathsf{LHS}$ is undecidable.
\end{theorem}

We show this by reduction of the $\mathbb{N}\times\mathbb{N}$ tiling problem \cite{tiling} to the satisfiability problem for $\mathsf{LHS}$. Let $\mathbb{T}= \{T_1,\dots,T_n\}$ be some fixed set of tile types. For each $T_i\in \mathbb{T}$, we use $\mathrm{up}(T_i)$, $\mathrm{down}(T_i)$, $\mathrm{left}(T_i)$ and $\mathrm{right}(T_i)$ to represent the colors of its up, down, left and right edges, respectively. We say that {\em $\mathbb{T}$ tiles $\mathbb{N}\times\mathbb{N}$} if there is a function $g:\mathbb{N}\times\mathbb{N}\to\mathbb{T}$ such that for all $n,m\in\mathbb{N}$,
\begin{center}
    $\mathrm{right}(g(n,m))=\mathrm{left}(g(n+1,m))$ and $\mathrm{up}(g(n,m))=\mathrm{left}(g(n,m+1))$.
\end{center}
Functions satisfying the conditions above are called {\em tiling functions}. In what follows, to show that $\LHS$ is undecidable, we present a formula $\phi_\mathbb{T}$ such that
\begin{center}
    $\phi_\mathbb{T}$ is satisfiable if and only if $\mathbb{T}$ tiles $\mathbb{N}\times\mathbb{N}$.
\end{center}
Let $\mathsf{Label}=\cset{u,r}\cup\cset{t_i:1\le i\le n}$ be a set of labels. Let $\mathsf{NV}^L=\cset{p^l:p\in\mathsf{Label}}$ and $\mathsf{NV}^R=\cset{p^r:p\in\mathsf{Label}}$ be sets of new variables. For convenience, we denote $\bigvee_{1\le i\le n}t^l_i$ by $t^l$ and $\bigvee_{1\le i\le n}t^r_i$ by $t^r$. We write $\D_u\phi$ for  $t^l\wedge\D(u^l\wedge\D(t^l\wedge\phi))$ and $\D_r\phi$ for $t^l\wedge\D(r^l\wedge\D(t^l\wedge\phi))$. Operators $\bd_u$ and $\bd_r$ are defined similarly. The dual of these operators are defined as usual, for example, $\B_u\phi:=\neg\D_u\neg\phi$.

The formula $\phi_\mathbb{T}$ is the conjunction of those in the groups below. To facilitate discussion, let $\M=(W,R,V)$ be a model and $w,v\in W$ s.t. $\M,w,v\md\phi_\mathbb{T}$.

\vspace{2mm}

\noindent \textbf{Group 1 (Basic requirements):}
\begin{enumerate}
    \item[] (SP) $I\wedge\B\B\bd I\wedge\D t^l$
\item[] (VL1) $\B\bb(I\to\bigwedge_{p\in\mathsf{Label}}(p^l\eq p^r))$
\item[] (VL2) $\B\bigwedge_{p\in\mathsf{Label}}(p^l\eq\bigwedge_{p\neq q\in\mathsf{Label}}\neg q^l)$
\end{enumerate}
Let us explain the meanings of the formulas in Group 1. Intuitively, we can treat $t,u,r$ as labels. The formula (SP) says that $w=v$, $R(R(w))\sub R(w)$ and there is some $v\in R(w)$ which is labelled by $t$.  (VL1) indicates that for any $s\in R(w)$, its  `left-label' and `right-label' are always the same. Moreover, (VL2) shows that every point $s\in R(w)$ has exactly one label.


\vspace{2mm}

\noindent
\textbf{Group 2 (Grid requirements):}
\begin{enumerate}
    \item[] (TU1) $\B\bb(t^l\wedge I\to\D (u^l\wedge\bb(u^r\to I)))$
    \item[] (TU2) $\B\bb(u^l\wedge I\to\D (t^l\wedge\bb(t^r\to I)))$
    \item[] (TR1) $\B\bb(t^l\wedge I\to\D (r^l\wedge\bb(r^r\to I)))$
    \item[] (TR2) $\B\bb(r^l\wedge I\to\D (t^l\wedge\bb(t^r\to I)))$
    \item[] (URT) $\B\bb(t^l\wedge I\to\B_u\bb_r\D_r\bd_u I)$
\end{enumerate}
We can assume that $\M$ is a model generated by $w\in W$ (Proposition \ref{prop:generated-model}). Let 
\begin{center}
    $R_u=\cset{\tup{s,t}\in R: \M,s,t\md t^l\wedge t^r\textit{ and for some } x\in V(u^l),\; sRx\textit{ and } xRt}$.
\end{center}
It follows from (TU1) and (TU2) that for all $s\in R(w)$, $|R_u(s)|=1$.\footnote{For any set $A$, we use $|A|$ for its {\em cardinality}.} Similarly, we can define $R_r$, and by (TR1) and (TR2),  for all $s\in R(w)$, $|R_r(s)|=1$. From (URT), we can infer that for all $v\in R(w)$, $R_r(R_u(v))=R_u(R_r(v))$.

\vspace{2mm}

\noindent
\textbf{Group 3 (Tiling the model):}
\begin{enumerate}
    \item[] (T1) $\B(t^l\to\bigwedge_{i=1}^n(t^l_i\to\D_u\bigvee_{1\le j\le n\;\&\; \mathrm{up}(T_i)=\mathrm{down}(T_j)}t^l_j))$;
    \item[] (T2) $\B(t^l\to\bigwedge_{i=1}^n(t^l_i\to\D_r\bigvee_{1\le j\le n\;\&\; \mathrm{right}(T_i)=\mathrm{left}(T_j)}t^l_j))$.
\end{enumerate}
The formulas in Group 3 are standard, which tell us that $\mathbb{T}$ `tiles' $R(w)\cap V(t^l)$.

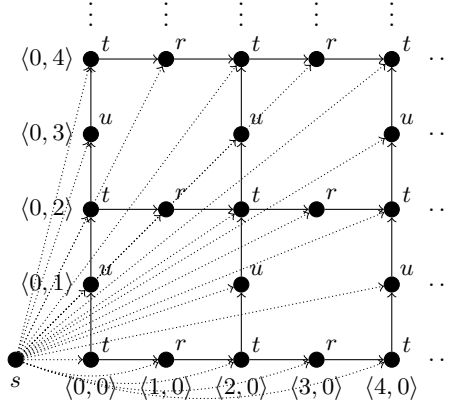
\begin{figure}
    \centering
    \begin{tikzpicture}
    \node(-10)[circle,draw,inner sep=0pt,minimum size=2mm,fill=black] at (-1,0)[label=below:$s$]{};
    \node(00)[circle,draw,inner sep=0pt,minimum size=2mm,fill=black] at (0,0)[label=below:$\tup{0,0}$]{};
    \node(10)[circle,draw,inner sep=0pt,minimum size=2mm,fill=black] at (1,0)[label=below:$\tup{1,0}$]{};
    \node(12)[circle,draw,inner sep=0pt,minimum size=2mm,fill=black] at (1,2){};
    \node(14)[circle,draw,inner sep=0pt,minimum size=2mm,fill=black] at (1,4){};
    \node(20)[circle,draw,inner sep=0pt,minimum size=2mm,fill=black] at (2,0)[label=below:$\tup{2,0}$]{};
    \node(21)[circle,draw,inner sep=0pt,minimum size=2mm,fill=black] at (2,1){};
    \node(22)[circle,draw,inner sep=0pt,minimum size=2mm,fill=black] at (2,2){};
    \node(23)[circle,draw,inner sep=0pt,minimum size=2mm,fill=black] at (2,3){};
    \node(24)[circle,draw,inner sep=0pt,minimum size=2mm,fill=black] at (2,4){};
    \node(30)[circle,draw,inner sep=0pt,minimum size=2mm,fill=black] at (3,0)[label=below:$\tup{3,0}$]{};
    \node(32)[circle,draw,inner sep=0pt,minimum size=2mm,fill=black] at (3,2){};
    \node(34)[circle,draw,inner sep=0pt,minimum size=2mm,fill=black] at (3,4){};
    \node(40)[circle,draw,inner sep=0pt,minimum size=2mm,fill=black] at (4,0)[label=below:$\tup{4,0}$]{};
    \node(41)[circle,draw,inner sep=0pt,minimum size=2mm,fill=black] at (4,1){};
    \node(42)[circle,draw,inner sep=0pt,minimum size=2mm,fill=black] at (4,2){};
    \node(43)[circle,draw,inner sep=0pt,minimum size=2mm,fill=black] at (4,3){};
    \node(44)[circle,draw,inner sep=0pt,minimum size=2mm,fill=black] at (4,4){};
    \node(50)[] at (4.7,0){$\dots$};
    \node(51)[] at (4.7,1){$\dots$};
    \node(52)[] at (4.7,2){$\dots$};
    \node(53)[] at (4.7,3){$\dots$};
    \node(54)[] at (4.7,4){$\dots$};
    \node(01)[circle,draw,inner sep=0pt,minimum size=2mm,fill=black] at (0,1)[label=left:$\tup{0,1}$]{};
    \node(02)[circle,draw,inner sep=0pt,minimum size=2mm,fill=black] at (0,2)[label=left:$\tup{0,2}$]{};
    \node(03)[circle,draw,inner sep=0pt,minimum size=2mm,fill=black] at (0,3)[label=left:$\tup{0,3}$]{};
    \node(04)[circle,draw,inner sep=0pt,minimum size=2mm,fill=black] at (0,4)[label=left:$\tup{0,4}$]{};
    \node(05)[] at (0,4.7){$\vdots$};
    \node(15)[] at (1,4.7){$\vdots$};
    \node(25)[] at (2,4.7){$\vdots$};
    \node(35)[] at (3,4.7){$\vdots$};
    \node(45)[] at (4,4.7){$\vdots$};
    \draw[->,densely dotted](-10) to (00);
    \draw[->,densely dotted](-10) to [bend right=20] (10);
    \draw[->,densely dotted](-10) to [bend right=20] (20);
    \draw[->,densely dotted](-10) to [bend right=20] (30);
    \draw[->,densely dotted](-10) to [bend right=20] (40);
    \draw[->,densely dotted](-10) to (01);
    \draw[->,densely dotted](-10) to (21);
    \draw[->,densely dotted](-10) to (41);
    \draw[->,densely dotted](-10) to (02);
    \draw[->,densely dotted](-10) to (12);
    \draw[->,densely dotted](-10) to (22);
    \draw[->,densely dotted](-10) to (32);
    \draw[->,densely dotted](-10) to (42);
    \draw[->,densely dotted](-10) to (03);
    \draw[->,densely dotted](-10) to (23);
    \draw[->,densely dotted](-10) to (43);
    \draw[->,densely dotted](-10) to (04);
    \draw[->,densely dotted](-10) to (14);
    \draw[->,densely dotted](-10) to (24);
    \draw[->,densely dotted](-10) to (34);
    \draw[->,densely dotted](-10) to (44);
    \draw[->](00) to (01);
    \draw[->](01) to (02);
    \draw[->](02) to (03);
    \draw[->](03) to (04);
    \draw[->](20) to (21);
    \draw[->](21) to (22);
    \draw[->](22) to (23);
    \draw[->](23) to (24);
    \draw[->](40) to (41);
    \draw[->](41) to (42);
    \draw[->](42) to (43);
    \draw[->](43) to (44);
    \draw[->](00) to (10);
    \draw[->](10) to (20);
    \draw[->](20) to (30);
    \draw[->](30) to (40);
    \draw[->](02) to (12);
    \draw[->](12) to (22);
    \draw[->](22) to (32);
    \draw[->](32) to (42);
    \draw[->](04) to (14);
    \draw[->](14) to (24);
    \draw[->](24) to (34);
    \draw[->](34) to (44);

    \draw (0.2,0.2) node{$t$};
    \draw (1.2,0.2) node{$r$};
    \draw (2.2,0.2) node{$t$};
    \draw (3.2,0.2) node{$r$};
    \draw (4.2,0.2) node{$t$};
    \draw (0.2,1.2) node{$u$};
    \draw (2.2,1.2) node{$u$};
    \draw (4.2,1.2) node{$u$};

    \draw (0.2,2.2) node{$t$};
    \draw (1.2,2.2) node{$r$};
    \draw (2.2,2.2) node{$t$};
    \draw (3.2,2.2) node{$r$};
    \draw (4.2,2.2) node{$t$};
    \draw (0.2,3.2) node{$u$};
    \draw (2.2,3.2) node{$u$};
    \draw (4.2,3.2) node{$u$};

    \draw (0.2,4.2) node{$t$};
    \draw (1.2,4.2) node{$r$};
    \draw (2.2,4.2) node{$t$};
    \draw (3.2,4.2) node{$r$};
    \draw (4.2,4.2) node{$t$};
    \end{tikzpicture}
    \caption{The model $\M_h$: Both dotted arrows and solid arrows represent the relation $R$.}
    \label{figure:spy}
\end{figure}

\begin{lemma}
    If $\mathbb{T}$ tiles $\mathbb{N}\times\mathbb{N}$, then $\phi_\mathbb{T}$ is satisfiable.
\end{lemma}
\begin{proof}
    Let $h:\mathbb{N}\times\mathbb{N}\to\mathbb{T}$ be a tiling function. Define $\M_h=(W,R,V)$ as follows:
    \begin{itemize}
        \item[$\bullet$] $W=W_0\cup\cset{s}$, where $W_0=\cset{\tup{n,m}\in\mathbb{N}\times\mathbb{N}:n\times m\textit{ is even}}$ 
        \item[$\bullet$]  $R=R_r\cup R_u\cup(\cset{s}\times W_0)$, where
        \begin{itemize}
            \item $R_r=\cset{\tup{\tup{k,2l},\tup{k+1,2l}}:k,l\in\mathbb{N}}$
            \item $R_u=\cset{\tup{\tup{2k,l},\tup{2k,l+1}}:k,l\in\mathbb{N}}$
        \end{itemize}
        \item[$\bullet$]  $V$ is a valuation such that
        \begin{itemize}
            \item[$\bullet$]  $V(r^l)=V(r^r)=\cset{\tup{2k+1,2l}\in W:k,l\in\mathbb{N}}$
            \item[$\bullet$] $V(u^l)=V(u^r)=\cset{\tup{2k,2l+1}\in W:k,l\in\mathbb{N}}$
            \item[$\bullet$]  $V(t^l_i)=V(t^r_i)=\cset{\tup{2k,2l}\in W:k,l\in\mathbb{N},h(k,l)=T_i}$ for all $1\le i\le n$.
            \item[$\bullet$]  $V(p^l)=V(q^r)=\emptyset$ for all other $p^l,q^r\in \mathsf{L}\cup \mathsf{R}$.
        \end{itemize}
    \end{itemize}
    The model $\M_h$ is shown in Fig \ref{figure:spy}. It is easy to verify that $\M_h,s,s\md\phi_\mathbb{T}$.\qed
\end{proof}

\begin{lemma}
    If $\phi_\mathbb{T}$ is satisfiable, then $\mathbb{T}$ tiles $\mathbb{N}\times\mathbb{N}$.
\end{lemma}
\begin{proof}
    Suppose $\M=(W,R,V)$ is a model generated by $s\in W$ and $\M,s,s\md\phi_\mathbb{T}$ (Proposition \ref{prop:generated-model}). It suffices to define a tiling function $g:\mathbb{N}\times\mathbb{N}\to T$. By (SP1), $V(t^l)\not=\emptyset$. Also, it follows from (TU1) and (TU2) that for each $w\in V(t^l)$, there is exactly one state $v\in V(t^l)$ such that $wRxRv$ for some $x\in V(u^l)$, and we denote the state $v$ by $\mathsf{up}(w)$. Then $\mathsf{up}:V(t^l)\to V(t^l)$ is a function. Similarly, due to (TR1) and (TR2), we can define a function $\mathsf{right}:V(t^l)\to V(t^l)$. Let $w_0\in V(t^l)$. We inductively define a function $g:\mathbb{N}\times\mathbb{N}\to V(t^l)$ as follows:
    \begin{align*}
        g(\tup{0,0})&=w_0,\\
        g(\tup{n,m+1})&=\mathsf{up}(g(\tup{n,m})),\\
        g(\tup{n+1,m})&=\mathsf{right}(g(\tup{n,m})).
    \end{align*}
   Now from  (URT) we can infer that  for each $w\in V(t^l)$, $\mathsf{up}(\mathsf{right}(w))=\mathsf{right}(\mathsf{up}(w))$. Then, for all $\tup{n,m}\in\mathbb{N}\times\mathbb{N}$, 
    {\small $$\mathsf{up}(g(\tup{n+1,m}))=\mathsf{up}(\mathsf{right}(g(\tup{n,m})))=\mathsf{right}(\mathsf{up}(g(\tup{n,m})))=\mathsf{right}(g(\tup{n,m+1})).$$}
    
 \noindent  Hence function $g$ is well-defined. Let $h:V(t^l)\to T$ be the function such that for each $1\le i\le n$, $h(w)=T_i$ if and only if $w\in V(t^l_i)$. Finally, by the formulas in Group 3, it is clear that $h\circ g:\mathbb{N}\times\mathbb{N}\to T$ is a tiling function. \qed
\end{proof}



\section{van Benthem Characterization Theorem}\label{sec:expressiveness}
This section is devoted to the expressive power of $\mathsf{LHS}$. Precisely, based on the notions of its first-order translation and bisimulations  developed in \cite{Li2021,LHS-journal}, we will provide a van Benthem style characterization theorem for the logic.

Let $\L^1$ be the first-order language consisting of the following: two disjoint countable sets $\mathsf{L}^1=\cset{P^l_i:i\in\mathbb{N}}$ and $\mathsf{R}^1=\cset{P^r_i:i\in\mathbb{N}}$ of unary predicates, a binary relation $R$ and equality $\equiv$. The first-order translation $\mathrm{T}_\tup{x,y}:\L\to\L^1$ for $\LHS$ is given in Appendix A. For any set $\Phi\subseteq \mathcal{L}$ of formulas, we define $\mathrm{T}_\tup{x,y}(\Phi):=\{\mathrm{T}_\tup{x,y}(\varphi):\varphi\in\Phi\}$. Moreover, notions of {\em $\LHS$-bisimulation} and {\em $\LHS$-saturation} are also introduced in Appendix A.

 Let $\Gamma(x_1,\dots,x_n)\subseteq\L^1$. We say that {\em $\M=(W,R,V)$ realizes $\Gamma$}  if there are $a_1,\dots,a_n\in W$ s.t. $\M\md\gamma[a_1,\dots,a_n]$ for all $\gamma\in\Gamma$. Also, let $A\sub W$. For each $a\in A$, let $c_a$ be a constant symbol. Let $\L^1_A=\L^1\cup\cset{c_a:a\in A}$ and let $\M_A$ denote the $\L^1_A$-expansion $\M$ such that for all $a\in A$, $a$ has the value $c_a$. 
 \begin{definition}
     A model $\M=(W,R,V)$ is {\em $\omega$-saturated}, if for all $A\subseteq W$, $\M_Y$ realizes every $\Gamma(x)\subseteq \L^1_Y$ whose finite subsets are all realized in $\M_Y$.
 \end{definition}

\begin{proposition}\label{prop:omegaSatu-LHSSatu}
All $\omega$-saturated models $\M=(W,R,V)$ are $\mathsf{LHS}$-saturated.
\end{proposition}
\begin{proof}
    Let $\Sigma\sub\L$ be finitely satisfiable in $R(w)\times\cset{v}$ and $w,v\in W$. Then let $\Delta(x)=\cset{Rc_wx}\cup\cset{\mathrm{T}_\tup{x,y}(\phi)[y/c_v]:\phi\in\Sigma}$. Every finite subset of $\Delta(x)$ is realized by some $u\in R(w)$ in $\M_\cset{w,v}$ (Proposition \ref{prop:ST} in Appendix A). Since $\M$ is $\omega$-saturated, $\Delta(x)$ is realized in $\M_\cset{w,v}$. So, there is some $u\in W$ such that $\M_\cset{w,v}\md\Delta(x)[u]$. Thus $u\in R(w)$ and $\tup{\M,u,v}\md\Sigma$.  Similarly, if $\Sigma$ is finitely satisfiable in $\cset{w}\times R(v)$, then $\Sigma$ is satisfiable in $\cset{w}\times R(v)$.  \qed
\end{proof}


\begin{corollary}\label{coro:OSandMEtoBis}
    Let $\M=(W,R,V)$ and $\M'=(W',R',V')$ be $\omega$-saturated models, $s,t\in W$ and $s',t'\in W'$. If $\tup{\M,s,t}$ and $\tup{\M',s',t'}$ satisfy the same $\mathcal{L}$-formulas, then $\tup{\M,s,t}\bis\tup{\M',s',t'}$.
\end{corollary}




Let $\M$ be a model, $\mathbb{I}$ a countable set and $U$ an incomplete ultrafilter over $\mathbb{I}$. Then we write $\prod_U\M$ for the ultrapower of $\M$ modulo $U$.\footnote{The notions {\em ultrafilter} and {\em ultrapower} are introduced in Appendix B.} 

\begin{proposition}\label{prop:FromMT}
    Let $\M=(W,R,V)$ be a model, $\mathbb{I}$ a countable set and $U$ an incomplete ultrafilter over $\mathbb{I}$. For each $w\in W$, let $f_w=\mathbb{I}\times\cset{w}$. Then,
    \begin{enumerate}
        \item $\prod_U\M$ is $\omega$-saturated.
        \item For any $\alpha(x,y)\in\L^1$ and $s,t\in W$,
                $\M\md\alpha[s,t]$ iff $\prod_U\M\md\alpha[(f_s)_U,(f_t)_U]$.
        \item For any $\L$-formula $\phi$ and $s,t\in W$,
            $\M,s,t\md\phi$ iff $\prod_U\M,(f_s)_U,(f_t)_U\md\phi$.
    \end{enumerate}
\end{proposition}
\begin{proof}
 The first item  follows from \cite[p.384, Theorem 6.1.1]{Chang1990_Chp4}. The second follows from \cite[p.217, Theorem 4.1.9]{Chang1990_Chp4}. The last one follows from the second item and Proposition \ref{prop:ST} in Appendix A immediately.\qed   
\end{proof}
 
We say that {\em an $\L^1$-formula $\alpha(x,y)$ is invariant for $\LHS$-bisimulation}, if for all $\tup{\M,s,t}$ and $\tup{\M',s',t'}$ that are $\LHS$-bisimilar, $\M\md\alpha[s,t]$ iff $\M'\md\alpha[s',t']$. Now we can provide a van Benthem style characterization theorem for $\LHS$:

\begin{theorem}
For any $\alpha(x,y)\in\L^1$,
        $\alpha(x,y)$ is invariant for $\LHS$-bisimulation if and only if $\md\alpha\eq\beta$ for some $\beta(x,y)\in\mathrm{T}_\tup{x,y}(\L)$.
\end{theorem}
\begin{proof}
    The right-to-left direction follows from Proposition \ref{proposition:bisimulation-to-equivalence} (see Appendix A). For the other direction, let $\alpha(x,y)\in\L^1$ be invariant for $\LHS$-bisimulation. Define  $\mathrm{modal}(\alpha):=\cset{\beta\in\mathrm{T}_\tup{x,y}(\L):\alpha\md\beta}.$
    We show $\mathrm{modal}(\alpha)\md\alpha$. Let $\M$ be a model such that $\M\md\mathrm{modal}(\alpha)[a,b]$. It suffices to show that $\M\md\alpha[a,b]$. Let $\Phi$ be a set of formulas defined by 
    \begin{center}
        $\Phi:=\cset{\beta(x,y)\in\mathrm{T}_\tup{x,y}(\L):\M\md\beta[a,b]}\cup\cset{\alpha(x,y)}$.
    \end{center} 
    We claim that $\Phi$ is satisfiable. Suppose $\Phi$ is not satisfiable. Then there is a finite $\Phi_0\sub\Phi$ with $\Phi_0\md\neg\alpha$, which entails $\alpha\md\neg\bigwedge\Phi_0$ and so $\M\md\neg\bigwedge\Phi_0[a,b]$. Note that $\bigwedge\Phi_0\in\Phi$, we see $\M\md\bigwedge\Phi_0[a,b]$, which is a contradiction. Thus, $\Phi$ is satisfiable and there is a model $\N$ and states $w,u$ with $\N\md\Phi[w,u]$. Then by Proposition \ref{prop:ST} in Appendix A, $\tup{\M,a,b}$ and $\tup{\N,w,u}$ satisfy the same $\mathcal{L}$-formulas. Let $U$ be an incomplete ultrafilter over $\mathbb{N}$. Then by Proposition \ref{prop:FromMT}(3), $\tup{\prod_U\M,(f_a)_U,(f_b)_U}$ and $\tup{\prod_U\N,(f_c)_U,(f_d)_U}$ satisfy the same $\mathcal{L}$-formulas. By Proposition \ref{prop:FromMT}(1) and Corollary \ref{coro:OSandMEtoBis}, $\tup{\prod_U\M,(f_a)_U,(f_b)_U}\bis\tup{\prod_U\N,(f_w)_U,(f_u)_U}$. Since $\N\md\alpha[w,u]$, by Proposition \ref{prop:FromMT}(2), $\prod_U\N\md\alpha[(f_w)_U,(f_u)_U]$. Since $\alpha(x,y)$ is invariant for $\LHS$-bisimulation, we have $\prod_U\M\md\alpha[(f_a)_U,(f_b)_U]$. By Proposition \ref{prop:FromMT}(2), $\M\md\alpha[a,b]$. Hence, $\mathrm{modal}(\alpha)\md\alpha$. By the Compactness Theorem, there is a finite $\Sigma\sub\mathrm{modal}(\phi)$ such that $\Sigma\md\alpha$. Then we see $\md\alpha\eq\bigwedge\Sigma$.\qed
\end{proof}

Finally, it is worthwhile to notice that when we restrict our attention to $\LHS^-$, by adapting  the arguments for $\LHS$, we can also obtain a characterization theorem for the expressiveness of $\LHS^-$, but we omit the details to save space.

\section{Axiomatization and Decidability of $\mathsf{LHS}^-$}\label{sec:axiomatization}

In this section, we turn our attention to $\mathsf{LHS}^-$. Precisely, we will provide a proof system for the logic, which is also helpful to show that its satisfiability problem is decidable. To achieve the former, instead of applying directly the usual techniques involving canonical models, we will make a detour: very roughly, we will first separate the `black part' and the `white part' of the language $\L^-$ of $\mathsf{LHS}^-$ and then build a desired calculus on that for the standard modal logic $\mathsf{K}$. The details will indicate that  containing two kinds of propositional variables in $\L^-$ makes $\mathsf{LHS}^-$ very different from its counterpart $\mathsf{K\times K}$ in product logic. Let us now introduce the details.

  A formula $\phi\in\L^-$ is called {\em a clean formula} if it contains only black modal operators or white modal operators. There are many formulas in the language $\L^-$ of $\mathsf{LHS}^-$ containing nested black modalities and white modalities. However, as we shall see, {\em every $\phi\in\L^-$ is logically equivalent to a Boolean combination of some clean formulas}.

\begin{definition}
    Languages $\L_\B$ and $\L_\bb$ are given by:
    \begin{align*}
        \L_\B\ni \phi &::= p^l \mid \neg\phi \mid \phi\wedge\phi \mid \B\phi,\\
        \L_\bb\ni \phi &::= p^r \mid \neg\phi \mid \phi\wedge\phi \mid \bb\phi,
    \end{align*}
        where $p^l\in \mathsf{L}$ and $p^r\in \mathsf{R}$.
\end{definition}

Let $\mathsf{K}_\B$ and $\mathsf{K}_\bb$ denote {\em the minimal modal logics with the languages $\L_\B$ and $\L_\bb$}, respectively. As the case for the standard modal logic $\mathsf{K}$, the satisfiability problems for both the logics are decidable (cf. \cite{open-mind}). Also, except the difference of the languages, their proof systems are the same as that for $\mathsf{K}$ \cite{BRV2001},  and we write ${\bf K}_{\B}$ and ${\bf K}_{\bb}$ for them. In what follows, we write $\md_1$ for {\em the usual one-dimensional satisfaction relation}. By induction on formulas, we can show that:
\begin{proposition}\label{prop:disjoint-union}
    Let $\M=(W,R,V)$ be a $\mathsf{K}_\B$-model and $\N=(U,S,V')$ a $\mathsf{K}_\bb$-model such that $W\cap U=\emptyset$. Let $\psi\in \L_\B$, $\gamma\in \L_\bb$. Then for all $s\in W$ and~$t\in U$,
    \begin{itemize}
        \item[(1)] $\M, s\md_1 \psi$ if and only if $\M\uplus\N,s,t\md \psi$, and 
        \item[(2)] $\N, t\md_1 \gamma$ if and only if $\M\uplus\N,s,t\md \gamma$,
    \end{itemize}
    where $\M\uplus\N$ is the $\mathsf{LHS}$-model defined by $\M\uplus\N=(W\cup U,R\cup S,V\cup V')$. The $\mathsf{LHS}$-model $\M\uplus\N$ is called the disjoint union of $\M$ and $\N$.
\end{proposition}

Let $\M$ be a $\mathsf{K}_\B$-model and $\N$ a $\mathsf{K}_\bb$-model. Since there are always isomorphic copies of them with disjoint domains, we can always assume that the domains of $\M$ and $\N$ are disjoint and construct the disjoint union of $\M$ and $\N$.

For an arbitrary $\LHS$-model $\M=(W,R,V)$, by restricting the valuation to $\mathsf{L}$ (we write $V|_{\mathsf{L}}$ for it), we can obtain a model $\M|_{\mathsf{L}}=(W,R,V|_{\mathsf{L}})$ for $\L_\B$, and similarly, by restricting $V$ to $\mathsf{R}$ (we write $V|_{\mathsf{R}}$ for it), we can get a model $\M|_{\mathsf{R}}=(W,R,V|_{\mathsf{R}})$ for $\L_\bb$. By induction on formulas, it is simple to prove that

\begin{proposition}\label{prop:models}
Let $\tup{\M,s,s}$ be a pointed $\mathsf{LHS}$-model. Then, for any $\phi\in\L_\B$,  $\M,s,s\md\phi$ iff $\M|_{\mathsf{L}},s\md_1\phi$. Also, for any $\phi\in\L_\bb$,  $\M,s,s\md\phi$ iff $\M|_{\mathsf{R}},s\md_1\phi$.
\end{proposition}

Before the next step, let us first recall some  concepts and facts about propositional logic. Let $\L_p$ denote the propositional language whose propositional variables come from $\mathsf{L}\cup\mathsf{R}$. For each $\phi\in\L_p$, we write $\phi(p_1,\dots,p_n)$ if the propositional variables occurring in $\phi$ are among $p_1,\dots,p_n$. Let $\phi(\alpha_1,\dots,\alpha_n)$ denote the formula obtained from $\phi(p_1,\dots,p_n)$ by simultaneously substituting $p_1,\dots,p_n$ with $\alpha_1,\dots,\alpha_n$ respectively. Let $\mathsf{PL}$ denote the set of all valid formulas in $\L_p$. A sound and complete Hilbert style calculus $\mathbf{PL}$ for $\mathsf{PL}$ can be given in a usual way.



\begin{definition}
A formula $\phi\in\L_p$ is {\em in conjunctive normal form (CNF)}, if $\phi$ is of the form
$\bigwedge_{i=1}^n\bigvee_{j=1}^{m_i}\phi_{ij}$,
where $n,m_1,\dots,m_n\in\mathbb{N}^+$ and each $\phi_{ij}$ is a propositional variable or a negation of a propositional variable. 
\end{definition}

We say that a formula $\phi$ is a {\em CNF-formula} if $\phi$ is in conjunctive normal form. Let $\mathsf{CNF}_p$ denote the set of all formulas $\phi\in\L_p$ in CNF.

\begin{proposition}\label{prop:CNF-prop}
    There is a  function $h:\L_p\to\mathsf{CNF}_p$ such that for all $\phi\in\L_p$, $\vdash_{\mathbf{PL}} \phi\eq h(\phi)$. 
\end{proposition}
\begin{proof}
    Such a function can be found in many textbooks of mathematical logic (see e.g., \cite[p. 221, Theorem 4.7]{mathmatical-logic}).
    \qed
\end{proof}

\begin{definition}\label{def:cleanFormula}
    Let $\phi\in\L^-$. Then we say that $\phi\in\L^-$ is a {\em clean formula}, if there are $\psi_1,\dots,\psi_n\in\L_\B$, $\gamma_1,\dots,\gamma_m\in\L_\bb$ and $\alpha(p^l_1, \dots, p^l_n, p^r_1,\dots, p^r_m)\in\L_p$ such that $\phi=\alpha(\psi_1,\dots,\psi_n,\gamma_1,\dots,\gamma_m)$.  Moreover, if $\alpha$ is in CNF, then $\phi$ is called a {\em clean CNF-formula}. Let $\L_c$ and $\mathsf{CNF}_c$ denote {\em the set of all clean formulas} and {\em the set of all clean CNF-formulas}, respectively.
\end{definition}


Table \ref{tab:proof-system-lhs^-} presents a Hilbert style calculus $\bf{LHS}^-$ for $\mathsf{LHS}^-$, which is a direct extension of the calculi $\bf{PL}$, ${\bf K}_{\B}$ and  ${\bf K}_{\bb}$. Therefore, we have the following: 

\begin{table}
    \centering
\begin{tabular}{|ll|}
\cline{1-2}
\multicolumn{2}{|c|}{ Proof  system $\bf{LHS}^-$ for $\mathsf{LHS}^-$}\\
\cline{1-2}
\multicolumn{2}{|l|}{\;\;Axioms:}\\
\cline{1-2}
\quad $(\mathrm{A1})$ &\quad  $p\to(q\to p)$,\; where $p,q\in \mathsf{L\cup R}$. \\
\quad $(\mathrm{A2})$ &\quad  $(p\to(q\to r))\to((p\to q)\to(p\to r))$,\; where $p,q,r\in \mathsf{L\cup R}$. \\
\quad $(\mathrm{A3})$ &\quad  $(\neg q\to\neg p)\to(p\to q)$, where $p,q\in \mathsf{L\cup R}$.   \\
\quad $(\mathrm{K}_{\B})$ &\quad  $\B(p^l\to q^l)\to(\B p^l\to\B q^l)$ \\
\quad $(\mathrm{K}_{\bb})$ &\quad  $\bb(p^r\to q^r)\to(\bb p^r\to\bb q^r)$ \\
\quad $(\mathrm{R}_{\B})$ &\quad  $\B(p^l\vee p^r)\eq(\B p^l\vee p^r)$ \\
\quad $(\mathrm{R}_{\bb})$ &\quad  $\bb(p^l\vee p^r)\eq( p^l\vee \bb p^r)$ \\
\cline{1-2}
\multicolumn{2}{|l|}{\;\;Inference rules:}\\
\cline{1-2}
\quad $(\mathrm{Sub})$ &\quad From $\alpha(p^l_1,\dots,p^l_n,p^r_1,\dots,p^r_m)$, infer $\alpha(\psi_1,\dots,\psi_n,\gamma_1,\dots,\gamma_m)$, \qquad \\
 & \quad where $\cset{\psi_i:1\leq i\leq n}\sub\L_\B$ and $\cset{\gamma_i:1\leq i\leq m}\sub\L_\bb$.\\
\quad $(\mathrm{MP})$ &\quad From $\alpha$ and $\alpha\to\beta$, infer $\beta$.\\
\quad $(\mathrm{Nec}_{\boxtimes})$ &\quad From $\alpha$, infer $\boxtimes\alpha$, for $\boxtimes\in\cset{\B,\bb}$. \\
\cline{1-2}
\end{tabular}    
    \caption{A proof  system $\bf{LHS}^-$ for $\mathsf{LHS}^-$}
    \label{tab:proof-system-lhs^-}
\end{table}

\begin{proposition}\label{prop:pc-to-hide-seek}
For all $\phi\in\L_p$, if $\vdash_{\bf{PL}}\phi$, then $\vdash_{\mathbf{LHS}^{-}}\phi$.
\end{proposition}

\begin{proposition}\label{prop:m-to-hide-seek}
For any formula $\varphi$ of $\L_\B$, if $\vdash_{\mathbf{K}_\B}\phi$, then $\vdash_{\mathbf{LHS}^{-}}\phi$.  Similarly, for any  formula $\varphi$ of  $\L_\bb$, if   $\vdash_{\mathbf{K}_\bb}\phi$, then $\vdash_{\mathbf{LHS}^{-}}\phi$.
\end{proposition}


\begin{corollary}\label{coro:Lc2CNFc}
    There is a function $f:\L_c\to\mathsf{CNF}_c$ such that for all $\phi\in\L_c$, it holds that $\vdash_{\mathbf{LHS}^{-}}\phi\eq f(\phi)$. Also, the resulting formula $f(\phi)$ is of the form $\bigwedge_{i=1}^n(\psi_i\vee\gamma_i)$, where $\bigwedge_{i=1}^n\psi_i\in\L_\B$ and $\bigwedge_{i=1}^n\gamma_i\in\L_\bb$.
\end{corollary}
\begin{proof}
    Let $\phi$ be a clean formula. Then, there are formulas $\beta(p_1,\dots,p_k)\in\L_p$ and $\alpha_1,\dots,\alpha_k\in\L_\B\cup\L_\bb$ such that $\phi=\beta(\alpha_1,\dots,\alpha_k)$. It follows from Proposition~\ref{prop:CNF-prop} that $\vdash_{\mathbf{PL}}\beta\eq h(\beta)$. Note that $h(\beta)$ is in CNF, and so $h(\beta)$ is of the form $\bigwedge_{i=1}^n\bigvee_{j=1}^{m_i}\beta_{ij}$ with $n,m_1,\dots,m_n\in\mathbb{N}^+$. For each $1\le i\le n$, we define:
    \begin{itemize}
        \item[$\bullet$]  $\psi_i:=(p^l\wedge\neg p^l)\vee\bigvee\cset{\beta_{ij}\in\L_\B:1\le j\le {m_i}}$, and
        \item[$\bullet$]  $\gamma_i:=(p^r\wedge\neg p^r)\vee\bigvee\cset{\beta_{ij}\in\L_\bb:1\le j\le {m_i}}$, 
    \end{itemize}
    where $p^l\in\mathsf{L}$ and $p^r\in\mathsf{R}$ are new propositional variables. Then it holds that
{\small{$$\vdash_{\mathbf{PL}}\bigwedge_{i=1}^n\bigvee_{j=1}^{m_i}\beta_{ij}\eq\bigwedge_{i=1}^n(\psi_i\vee\gamma_i).$$}}
By Proposition \ref{prop:pc-to-hide-seek}, 
    $\vdash_{\mathbf{LHS}^{-}}\beta(p_1,\dots,p_k)\eq\bigwedge_{i=1}^n(\psi_i\vee\gamma_i)$.
    Next, applying the inference rule (Sub) of the calculus $\mathbf{LHS}^{-}$, we can obtain
{\small{$$\vdash_{\mathbf{LHS}^{-}}\phi\eq\bigwedge_{i=1}^n(\psi_i(\alpha_1,\dots,\alpha_k,p^l)\vee\gamma_i(\alpha_1,\dots,\alpha_k,p^r)).$$}}

    Now, we can define a desired function $f:\L_c\to\mathsf{CNF}_c$ in the following manner:
    {\small{$$f(\phi)=\bigwedge_{i=1}^n(\psi_i(\alpha_1,\dots,\alpha_k,p^l)\vee\gamma_i(\alpha_1,\dots,\alpha_k,p^r)),$$}} which completes the proof.\qed
\end{proof}

\begin{lemma}\label{lem:clean-sat}
    Let $\M=(W,R,V)$ be a model and $s,t\in W$. Then, for all formulas $\phi\in\L_\B$ and $\psi\in\L_\bb$, the following equivalences hold:
    \begin{itemize}
        \item[(1)] $\M,s,t\md\phi$ if, and only if, for all $t'\in W$, $\M,s,t'\md\phi$.
        \item[(2)] $\M,s,t\md\psi$ if, and only if,  for all $s'\in W$, $\M,s',t\md\phi$.
        \item[(3)] $\M,s,t\md\B(\phi\vee\psi)$ if, and only if, $\M,s,t\md\B\phi\vee\psi$.
        \item[(4)] $\M,s,t\md\bb(\phi\vee\psi)$ if, and only if, $\M,s,t\md\phi\vee\bb\psi$.
    \end{itemize}
\end{lemma}
The proof is given in Appendix C. From the items (3) and (4) of  Lemma \ref{lem:clean-sat}, it is a matter of direct checking that:

\begin{lemma}\label{lem:Box-vee-distribute}
    For all $\phi\in\L_\B$ and $\psi\in\L_\bb$, both the formulas  $\B(\phi\vee\psi)\eq(\B\phi\vee\psi)$ and $\bb(\phi\vee\psi)\eq(\phi\vee\bb\psi)$ are valid. 
\end{lemma}


Now we move to showing the soundness of the proof system ${\bf LHS}^{-}$:

\begin{theorem}\label{theorem:soundness}
    For any formula $\phi\in\L^-$, $\vdash_{{\bf LHS}^{-}}\phi$ implies $\md\phi$.
\end{theorem}

\begin{proof}
  The validity of (A1), (A2), (A3), $(\mathrm{K}_\B)$ and $(\mathrm{K}_\bb)$ is easy to see. Also, Lemma \ref{lem:Box-vee-distribute} indicates that $(\mathrm{R}_\B)$ and $(\mathrm{R}_\bb)$ are valid. Moreover,  all inference rules of ${\bf LHS}^{-}$ preserve validity, and the details are left as an exercise. \qed
\end{proof}

Next, we consider for the completeness of the calculus ${\bf LHS}^{-}$.

\begin{definition}\label{def:cleanCNF-companion}
For any $\L^-$-formula $\varphi$, we define its {\em clean CNF companion} $\phi_c$ in the following inductive manner: 
{ \small{\begin{align*}
    (p^l)_c&:=p^l\vee(p^r_0\wedge\neg p^r_0), \text{ where } p^r \text{ is a new propositional variable from } \mathsf{R}. \\
        (p^r)_c&:=(p^l_0\wedge\neg p^l_0)\vee p^r, \text{ where } p^l \text{ is a new propositional variable from } \mathsf{L}.\\
        (\neg\phi)_c&:=f(\neg\phi_c)\\
        (\phi\wedge\psi)_c&:=\phi_c\wedge\psi_c\\
        (\B\phi)_c&:=\bigwedge_{i=1}^n(\B\psi_i\vee\gamma_i), \text{ where } \bigwedge_{i=1}^n\psi_i\in\L_\B, \bigwedge_{i=1}^n\gamma_i\in\L_\bb \text{ and } \phi_c=\bigwedge_{i=1}^n(\psi_i\vee\gamma_i).\\
        (\bb\phi)_c&:=\bigwedge_{i=1}^n(\psi_i\vee\bb\gamma_i), \text{ where } \bigwedge_{i=1}^n\psi_i\in\L_\B, \bigwedge_{i=1}^n\gamma_i\in\L_\bb \text{ and } \phi_c=\bigwedge_{i=1}^n(\psi_i\vee\gamma_i).
    \end{align*}}}
\end{definition}

\begin{example} Let us consider an example $\bb p^l$. With the clauses above, we have
    $(\bb p^l)_c=p^l\vee\bb(p^r_0\wedge\neg p^r_0)$. Note that if we define $(p^l)_c$ to be $p^l$, then we cannot ensure that a formula  and its clean companion are equivalent: for instance, one can easily find a model making $\bb p^l \eq p^l$ false. Similarly for the clause of $p^r\in\mathsf{R}$.   
\end{example}

\begin{theorem}\label{thm:CNF}
  Let $\phi$ be a formula of $\L^-$. Then, its clean CNF companion $\phi_c$ is of the form $\bigwedge_{i=1}^n(\psi_i\vee\gamma_i)$ with $\bigwedge_{i=1}^n\psi_i\in\L_\B$ and $\bigwedge_{i=1}^n\gamma_i\in\L_\bb$. Moreover, it holds that $\vdash_{{\bf LHS}^{-}}\phi\eq\phi_c$.
\end{theorem}
\begin{proof}
The first part of the theorem is easy, since $\phi_c$ is  a clean CNF formula. In what follow, by induction on $\phi\in\L^{-}$, we will show that $\vdash_{{\bf LHS}^{-}}\phi\eq\phi_c$. The cases for propositional atoms and $\wedge$ are straightforward, and we consider others.

  \vspace{2mm}

(1). First,  we consider $\phi=\neg\psi$. By the induction hypothesis, it holds that $\vdash_{{\bf LHS}^{-}}\psi\eq\psi_c$. So,  $\vdash_{{\bf LHS}^{-}}\neg\psi\eq\neg\psi_c$. Clearly, $\neg\psi_c\in \L_c$. From Corollary \ref{coro:Lc2CNFc} we know that $\vdash_{{\bf LHS}^{-}}  f(\neg\psi_c)\eq\neg\psi_c$. Also, with the clause for $\neg$ in Definition \ref{def:cleanCNF-companion}, we have $(\phi)_c=f(\neg\psi_c)$. Immediately, we obtain $\vdash_{{\bf LHS}^{-}}\phi_c\eq\phi$, as desired. 

 \vspace{2mm}

(2). Next, we consider $\phi=\B\psi$. Assume that $\psi_c=\bigwedge_{i=1}^n(\psi_i'\vee\gamma_i')$, where $\bigwedge_{i=1}^n\psi_i'\in\L_\B$ and $\bigwedge_{i=1}^n\gamma_i'\in\L_\bb$. Then,  $\vdash_{{\bf LHS}^{-}}\B\psi_c\eq\bigwedge_{i=1}^n\B(\psi_i'\vee\gamma_i').$ For simplicity, we write $\gamma(p^l,p^r)$ for $\B(p^l\vee p^r)\eq(\B p^l\vee p^r)$, which is exactly the axiom (R$_{\B}$). Note that for each $1\le i\le n$, $\psi_i'\in\L_\B$ and $\gamma_i'\in\L_\bb$. So, for each $1\le i\le n$, using the inference rule (Sub), we can obtain 
$\vdash_{{\bf LHS}^{-}}\gamma(\psi_i',\gamma_i')$, i.e., $\vdash_{{\bf LHS}^{-}}\B(\psi_i'\vee\gamma_i')\eq(\B\psi_i'\vee\gamma_i')$. Therefore, $\vdash_{{\bf LHS}^{-}}\bigwedge_{i=1}^n\B(\psi_i'\vee\gamma_i')\eq\bigwedge_{i=1}^n(\B\psi_i'\vee\gamma_i')$, which entails $\vdash_{{\bf LHS}^{-}}\B\psi_c\eq\phi_c$. By induction hypothesis, $\vdash_{{\bf LHS}^{-}}\psi\eq\psi_c$ and so $\vdash_{{\bf LHS}^{-}}\phi\eq\B\psi_c$. Hence $\vdash_{{\bf LHS}^{-}}\phi\eq\phi_c$. 

  \vspace{2mm}

 (3). Finally, the case for $\phi=\bb\psi$ is similar to (2). The proof is completed.
\qed
\end{proof}

\begin{example}
Applications of Theorem  \ref{thm:CNF} can be diverse. As an example, we show how to use it to prove $\vdash_{{\bf LHS}^{-}}\B\bb\neg\phi\eq\bb\B\neg\phi$. Let $(\neg\phi)_c=\bigwedge_{i=1}^n(\psi_i\vee\gamma_i)$ with $\bigwedge_{i=1}^n\psi_i\in\L_\B$ and $\bigwedge_{i=1}^n\gamma_i\in\L_\bb$. By Theorem \ref{thm:CNF}, $\vdash_{{\bf LHS}^{-}}\neg\phi\eq(\neg\phi)_c$. Then, it is simple to see that $\vdash_{{\bf LHS}^{-}}\bb\B\neg\phi\eq\bb\B(\neg\phi)_c$ and $\vdash_{{\bf LHS}^{-}}\B\bb\neg\phi\eq\B\bb(\neg\phi)_c$. Also, we can show  both
    $\vdash_{{\bf LHS}^{-}}\bb\B(\neg\phi)_c\eq\bigwedge_{i=1}^n(\B\psi_i\vee\bb\gamma_i)$ and 
    $\vdash_{{\bf LHS}^{-}}\B\bb(\neg\phi)_c\eq\bigwedge_{i=1}^n(\B\psi_i\vee\bb\gamma_i)$.
Thus, we obtain $\vdash_{{\bf LHS}^{-}}\B\bb\neg\phi\eq\bb\B\neg\phi$.
\end{example}

Now, with the help of Theorem \ref{thm:CNF}, we can show that ${\bf LHS}^{-}$ is complete with respect to finite $\mathsf{LHS}$-models:

\begin{theorem}\label{thm:completeness}
    For each $\phi\in\L^-$, if $\mathsf{Mod}^{<\omega}\md\phi$ then $\vdash_{{\bf LHS}^{-}}\phi$, where  $\mathsf{Mod}^{<\omega}$ is the class of  finite models. 
\end{theorem}

    \begin{proof}
        Let $\phi\in\L^-$ and $\not\vdash_{\mathbf{LHS}^-}\phi$. By Theorem \ref{thm:CNF}, $\not\vdash_{\mathbf{LHS}^-}\phi_c$ and $\phi_c$ is of the form $\bigwedge_{i=1}^n(\psi_i\vee\gamma_i)$ where $\bigwedge_{i=1}^n\psi_i\in\L_\B$ and $\bigwedge_{i=1}^n\gamma_i\in\L_\bb$. Since $\not\vdash_{\mathbf{LHS}^-}\phi_c$, there is some $i$ such that $\not\vdash_{\mathbf{LHS}^{-}}\psi_i\vee\gamma_i$. By Proposition~\ref{prop:m-to-hide-seek}, we have $\not\vdash_{\mathbf{K}_\B}\psi_i$ and $\not\vdash_{\mathbf{K}_\bb}\gamma_i$. By the completeness of $\mathbf{K}_\B$ and $\mathbf{K}_\bb$, both $\neg\psi_i$ and $\neg\gamma_i$ is satisfiable. Since $\mathbf{K}_\B$ and $\mathbf{K}_\bb$ have the finite model property, there are finite pointed $\mathbf{K}_\B$-model $\tup{\M,s}$ and finite pointed $\mathbf{K}_\bb$-model $\tup{\N,t}$ such that $\M,s\not\md\psi_i$ and $\N,t\not\md\gamma_i$. Then, by Proposition \ref{prop:disjoint-union}, $\M\uplus\N,s,t\md\neg\psi_i\wedge\neg\gamma_i$, which entails $\M\uplus\N,s,t\md\neg\phi_c$. By Theorem \ref{theorem:soundness} and Theorem \ref{thm:CNF}, $\mathsf{Mod}^{<\omega}\not\md\phi$. \qed
    \end{proof}



\begin{theorem}\label{thm:fmp}
    $\mathsf{LHS}^-$ enjoys the finite model property, and it is decidable.
\end{theorem}




\section{Conclusion}\label{sec:conclusion}
{\em Summary}\; The article is a technical continuation of \cite{Li2021,LHS-journal}, which explored the properties $\LHS$, a tool  to reason about the games of hide and seek. In the paper, we show that the satisfiability problem for $\LHS$ with a single relation is undecidable. Also, based on existing notions of bisimulations and first-order translation for $\LHS$, we provide a van Benthem style characterization theorem for the logic. Moreover, we develop a Hilbert style calculus for $\LHS^-$ and prove that its satisfiability problem  is decidable, and the details of our proofs are helpful to understand the differences between $\LHS^-$ and $\mathsf{K\times K}$. All these results can be transferred to logics generalizing $\LHS$ and $\LHS^-$ for games with $n>2$ players. 

\vspace{2mm}

\noindent {\em Related Work}\; As stated, our work is closely related to product logics $\mathsf{K\times K}$ \cite{many-dimensional-book} and $\mathsf{K\times^{\delta} K}$ \cite{product-undecidable,Kikot2010,diagonal-fmp}, cylindric modal logics \cite{Yde-thesis} and cylindric algebra \cite{cylindric-algebra}. Also, there is a line of logical investigation for the hide and seek game. Needless to say, the most relevant ones are \cite{Li2021} and its extension \cite{LHS-journal}. The latter offers a notion of bisimulation for $\LHS$ and its first-order translation, proves the undecidability of $\LHS$ with multiple relations, shows that the model-checking problems for both $\LHS$ and $\LHS^-$ are $\mathsf{P}$-complete, and identifies the counterpart of $\LHS$ in product logic. Moreover, \cite{hlhs-axiomatization} extends $\LHS$ and $\LHS^-$ with components from hybrid logics, and study the expressiveness and axiomatization of the resulting logics. Also, \cite{LHS-epistemic} develops logical tools to capture how players update their knowledge about other players' positions in the hide and seek game, in which players have only imperfect information. Finally, it is important to notice that besides the efforts made for the hide and seek game, many other graph games and their matching modal logics have also been studied in recent years, and we refer to \cite{graph-game-logic-design} for a broad research program on this topic and refer to \cite{Dazhu-thesis}  for extensive references to modal logics for graph games.

\vspace{2mm}

\noindent{\em Further Directions}\;  Except what we have explored in the article, there are a number of directions deserving  to be explored in future. A natural next step is to study the axiomatizability of $\LHS$, for which it might be useful to analyze the techniques developed for $\mathsf{K\times^{\delta} K}$ \cite{Kikot2010}. Closely related to this, \cite{Li2021,LHS-journal} provide preliminary discussions on the difference between the frameworks of $\LHS$ and $\LHS^-$ and product logics, but their exact relation remains to be explored. Also, besides the expressive power of $\LHS$ w.r.t. models, the equality constant $I$ improves the frame definability of $\LHS$ as well,\footnote{For instance, although {\em confluence} (i.e., $\forall x\forall y_1\forall y_2(Rxy_1\land Rxy_2\to\exists z(Ry_1z\land Ry_2z))$) is not definable with the basic modal language, using the techniques of {\em frame correspondence} \cite{correspondence-theory}, one can check that the  property can be simply defined as $I\to \B\bb\D \bd I$.} and it remains to have a comprehensive understanding of the expressive power of $\LHS$ w.r.t. frames. Moreover, another direction is to develop the proof theory for our logics, and for instance, provide sequent calculi and tableau systems for them. Finally, it is worthwhile to generalize our results to broader settings and consider the extensions of $\LHS$ and $\LHS^-$ with further operators, such as {\em graded modalities} \cite{graded-modalities,graded-modalities-Katsuhiko} to talk about the degree of a state in a graph.

\bibliographystyle{splncs04}
\bibliography{my}

\begin{thebibliography}{10}
\providecommand{\url}[1]{\texttt{#1}}
\providecommand{\urlprefix}{URL }
\providecommand{\doi}[1]{https://doi.org/#1}

\bibitem{open-mind}
van Benthem, J.: Modal Logic for Open Minds. CSLI Publications (2010)

\bibitem{graph-game-logic-design}
van Benthem, J., Liu, F.: Graph games and logic design. In: Liu, F., Ono, H.,
  Yu, J. (eds.) Knowledge, Proof and Dynamics, p. 125–146. Springer,
  Singapore (2020)

\bibitem{correspondence-theory}
van Benthem, J.: Correspondence theory. In: Gabbay, D., Guenthner, F. (eds.)
  Handbook of Philosophical Logic, Synthese Library (Studies in Epistemology,
  Logic, Methodology, and Philosophy of Science), vol.~165, pp. 167--247.
  Springer (1984)

\bibitem{BRV2001}
Blackburn, P., Rijke, M.d., Venema, Y.: Modal Logic. Cambridge University
  Press, Cambridge (2001)

\bibitem{Chang1990_Chp4}
Chang, C., Keisler, H.: Model Theory, Studies in Logic and the Foundations of
  Mathematics, vol.~73. Elsevier, Amsterdam, second edn. (1990)

\bibitem{mathmatical-logic}
Ebbinghaus, H.D., Flum, J., Thomas, W.: Mathematical Logic. Springer, Cham,
  third edn. (2021)

\bibitem{graded-modalities}
Fine, K.: In so many possible worlds. Notre Dame Journal of Formal Logic
  \textbf{13},  516--520 (1972)

\bibitem{many-dimensional-book}
Gabbay, D., Kurucz, A., Wolter, F., Zakharyaschev, M.: Many-Dimensional Modal
  Logics: Theory and Applications. Elsevier, Amsterdam (2003)

\bibitem{product-undecidable}
Hampson, C., Kikot, S., Kurucz, A.: The decision problem of modal product
  logics with a diagonal, and faulty counter machines. Studia Logica
  \textbf{104},  455--486 (2016)

\bibitem{cylindric-algebra}
Henkin, L., Monk, J.D., Tarski, A.: Cylindric Algebras, Part I, Studies in
  Logic and the Foundations of Mathematics, vol.~64. North-Holland, Amsterdam
  (1971)

\bibitem{Kikot2010}
Kikot, S.: Axiomatization of modal logic squares with distinguished diagonal.
  Mathematical Notes  \textbf{88},  238--250 (2010)

\bibitem{diagonal-fmp}
Kurucz, A.: Products of modal logics with diagonal constant lacking the finite
  model property. In: Ghilardi, S., Sebastiani, R. (eds.) Frontiers of
  Combining Systems. LNCS, vol.~5749, pp. 279--286 (2009)

\bibitem{Dazhu-thesis}
Li, D.: Formal Threads in the Social Fabric: Studies in the Logical Dynamics of
  Multi-Agent Interaction. Ph.D. thesis, Department of Philosophy, Tsinghua
  University and ILLC, University of Amsterdam (2021)

\bibitem{LHS-epistemic}
Li, D., Ghosh, S., Liu, F.: Dynamic-epistemic logic for the cops and robber
  game. Manuscript (2023)

\bibitem{Li2021}
Li, D., Ghosh, S., Liu, F., Tu, Y.: On the subtle nature of a simple logic of
  the hide and seek game. In: Silva, A., Wassermann, R., de~Queiroz, R. (eds.)
  Logic, Language, Information, and Computation. pp. 201--218. Springer, Cham
  (2021)

\bibitem{LHS-journal}
Li, D., Ghosh, S., Liu, F., Tu, Y.: A simple logic of the hide and seek game.
  Studia Logica  (2023). \doi{10.1007/s11225-023-10039-4}

\bibitem{cop-robber}
Nowakowski, R., Winkler, R.P.: Vertex-to-vertex pursuit in a graph. Discrete
  Math  \textbf{13},  235--239 (1983)

\bibitem{tiling}
Robinsonn, R.: Undecidability and nonperiodicity for tilings of the plane.
  Inventiones Mathematicae  \textbf{12},  177–209 (1971)

\bibitem{hlhs-axiomatization}
Sano, K., Liu, F., Li, D.: Hybrid logic of the hide and seek game. Manuscript
  (2023)

\bibitem{graded-modalities-Katsuhiko}
Sano, K., Ma, M.: Goldblatt-{T}homason-style theorems for graded modal
  language. In: Beklemishev, L., Goranko, V., Shehtman, V. (eds.) Advances in
  Modal Logic, vol.~8, pp. 330--349. CSLI Publications, San Diego (2010)

\bibitem{Yde-thesis}
Venema, Y.: Many-Dimensional Modal Logics. Ph.D. thesis, Universiteit van
  Amsterdam (1991)

\end{thebibliography}

\appendix

\section*{Appendix A: FO-Translation and Bisimulations for $\LHS$}

\begin{definition}[\cite{LHS-journal}]
    For any two variables $x$ and $y$, we  define {\em the first-order translation $\mathrm{T}_\tup{x,y}:\L\to\L^1$} recursively as follows:
\begin{center}
       $\mathrm{T}_\tup{x,y}(p^l_i):= P^l_ix$\qquad 
           $ \mathrm{T}_\tup{x,y}(p^r_i):=P^r_iy$\qquad
           $ \mathrm{T}_\tup{x,y}(I):={(x\equiv y)}$\vspace{1.5mm}\\ 
         $\mathrm{T}_\tup{x,y}(\neg\phi):=\neg\mathrm{T}_\tup{x,y}(\phi)$\qquad
            $\mathrm{T}_\tup{x,y}(\phi\wedge\psi):=\mathrm{T}_\tup{x,y}(\phi)\wedge\mathrm{T}_\tup{x,y}(\psi)$\vspace{1.5mm}\\
           $  \mathrm{T}_\tup{x,y}(\B\phi):=\forall{z}(Rxz\to\mathrm{T}_\tup{z,y}(\phi))$ \quad
           $ \mathrm{T}_\tup{x,y}(\bb\phi):=\forall{z}(Ryz\to\mathrm{T}_\tup{x,z}(\phi))    $   
\end{center}    
\end{definition}

 Now, the following result indicates the correctness of the translation:

    \begin{proposition}[\cite{LHS-journal}]\label{prop:ST}
        For any pointed $\LHS$-model $\tup{\M,s,t}$ and $\phi\in\L$,
        \begin{center}
            $\M,s,t\md\phi$ if and only if $\M\md\mathrm{T}_\tup{x,y}(\phi)[s,t]$.\footnote{By $\M\md\mathrm{T}_\tup{x,y}(\phi)[s,t]$, we mean that when values of $x,y$ in  $\mathrm{T}_\tup{x,y}(\phi)$ are $s,t$ respectively, $\mathrm{T}_\tup{x,y}(\phi)$ is satisfied by $\M$.}
        \end{center}
    \end{proposition}

    \begin{definition}[\cite{Li2021,LHS-journal}]
    Let $\M=(W,R,V)$ and $\M'=(W',R',V')$ be  models. A binary relation $Z\sub (W\times W)\times(W'\times W')$ is called {\em an $\mathsf{LHS}$-bisimulation between $\M$ and $\M'$}, notation $Z:\M\bis\M'$, if the following conditions hold for all $s,t,v\in W$ and $s',t',v'\in W'$:
    \begin{itemize}
        \item[$\bullet$] If $\tup{s,t}Z\tup{s',t'}$, then for all $p\in\mathsf{L}\cup\mathsf{R}$, $\M,s,t\md p$ if and only if $\M',s',t'\md p$.
        \item[$\bullet$] If $\tup{s,t}Z\tup{s',t'}$ and $v\in R(s)$, then there is $v'\in R'(s')$ s.t. $\tup{v,t}Z\tup{v',t'}$.
        \item[$\bullet$] If $\tup{s,t}Z\tup{s',t'}$ and $v\in R(t)$, then there is $v'\in R'(t')$ s.t. $\tup{s,v}Z\tup{s',v'}$.
        \item[$\bullet$] If $\tup{s,t}Z\tup{s',t'}$ and $v'\in R'(s')$, then there is $v\in R(s)$ s.t. $\tup{v,t}Z\tup{v',t'}$.
        \item[$\bullet$] If $\tup{s,t}Z\tup{s',t'}$ and $v'\in R'(t')$, then there is $v\in R(t)$ s.t. $\tup{s,v}Z\tup{s',v'}$.
        \item[$\bullet$] If $\tup{s,t}Z\tup{s',t'}$, then $s=t$ if and only if $s'=t'$.
    \end{itemize}
\end{definition}

If there is an $\mathsf{LHS}$-bisimulation $Z:\M\bis\M'$ s.t. $\tup{s,t}Z\tup{s',t'}$, then we say that {\em $\tup{\M,s,t}$ is $\mathsf{LHS}$-bisimular to $\tup{\M',s',t'}$} and write $\tup{\M,s,t}\bis\tup{\M',s',t'}$. 

\begin{proposition}[\cite{Li2021,LHS-journal}]\label{proposition:bisimulation-to-equivalence}
If $\tup{\M, s, t}\bis\tup{\M', s', t'}$, then $\tup{\M, s, t}$ and $\tup{\M', s', t'}$ satisfy the same $\LHS$-formulas. 
\end{proposition}

The converse of Proposition \ref{proposition:bisimulation-to-equivalence} holds for {\em $\LHS$-saturated models}:

\begin{definition}[\cite{Li2021,LHS-journal}]
    Let $\M=(W,R,V)$ be a model. A set $\Delta$ of formulas is said to be satisfiable in $X\sub W\times W$ if $\M,s,t\md\Delta$ for some $\tup{s,t}\in X$. Then $\M$ is said to be {\em $\mathsf{LHS}$-saturated} if for all $\Phi\sub\L$ and $w,v\in W$:
    \begin{enumerate}
        \item  If every finite subset $\Sigma$ of $\Phi$ is satisfiable in $R(w)\times\cset{v}$, then $\Phi$ is satisfiable in $R(w)\times\cset{v}$.
        \item If every finite subset $\Sigma$ of $\Phi$ is satisfiable in $\cset{w}\times R(v)$, then $\Phi$ is satisfiable in $\cset{w}\times R(v)$.
    \end{enumerate}
\end{definition}

\begin{proposition}[\cite{Li2021,LHS-journal}]\label{proposition:equivalenc-to-bisimulation}
For all $\M$ and $\M'$ that are $\LHS$-saturated, if $\tup{\M,s,t}$ and $\tup{\M',s',t'}$ satisfy the same formulas of $\LHS$, then $\tup{\M,s,t}\bis\tup{\M',s',t'}$. 
\end{proposition}

\section*{Appendix B: Ultrafilter, Ultraproduct, Ultrapower}


\begin{definition}
$U\in\mathcal{P}(\mathcal{P}(W))$ is  an {\em ultrafilter} over a non-empty set $W$ if, 
        \begin{itemize}
            \item[$\bullet$] $W\in U$;
            \item[$\bullet$] $X,Y\in U$ implies $X\cap Y\in U$;
            \item[$\bullet$] if $X\in U$ and $X\sub Y\sub W$, then $Y\in U$;
            \item[$\bullet$] for all $X\sub W$, $X\in U$ if and only if $(W\sm X)\not\in U$.
        \end{itemize}
       An ultrafilter is {\em incomplete} if it is not closed under countable intersections.
    \end{definition}
    
    \begin{definition}
        Let $\cset{\M_i=(W_i,R_i,V_i):i\in \mathbb{I}}$ be a family of models and $U$ an ultrafilter over $\mathbb{I}$. Then we write $\prod_{i\in\mathbb{I}}W_i$ for the set of all functions $f:\mathbb{I}\to\bigcup_{i\in\mathbb{I}} W_i$ such that $f(i)\in W_i$ for all $i\in\mathbb{I}$. For each $f\in\prod_{i\in\mathbb{I}}W_i$, we define 
        \begin{center}
            $f_U=\cset{g\in\prod_{i\in\mathbb{I}}W_i:\cset{i\in \mathbb{I}:f(i)=g(i)}\in U}$.
        \end{center}
        Then the {\em ultraproduct $\prod_U\M_i=(W_U,R_U,V_U)$ of $(\M_i)_{i\in\mathbb{I}}$ modulo $U$} is the model given by the following:
        \begin{enumerate}
            \item $W_U=\prod_UW_i=\cset{f_U:f\in\prod_{i\in \mathbb{I}}W_i}$.
            \item For each $f_U\in W_U$ and $p\in\mathsf{L}\cup\mathsf{R}$, $f_U\in V_U(p)$ iff $\cset{i\in \mathbb{I}:f(i)\in V_i(p)}\in U$.
            \item For all $f_U, g_U\in W_U$, $R_Uf_Ug_U$ iff $\cset{i\in \mathbb{I}:R_if(i)g(i)}\in U$.
        \end{enumerate}
        If $\M_i=\M$ for all $i\in \mathbb{I}$, then $\prod_U\M$ is called the {\em ultrapower of $\M$ modulo $U$}.
    \end{definition}

    \section*{Appendix C: Proof for Lemma \ref{lem:clean-sat}}

    \begin{proof}
        The proofs for items (1) and (2) are by induction on the complexity of $\phi$ and $\psi$, respectively. We omit the details for them. In what follows, we merely consider for (3), since (4) can be prove in a similar way.  

        For the direction from left to right, we prove the contrapositive and assume that  $\M,s,t\not\md\B\phi\vee\psi$. Then, $\M,s,t\md\neg\psi$ and there is some state $s'\in R(s)$ such that $\M,s',t\not\md\phi$. Note that $\neg\psi\in\L_\bb$. Now, using the item (2), we can obtain $\M,s',t\md\neg\psi$. Thus, it holds that $\M,s',t\not\md\phi\vee\psi$, and so $\M,s,t\not\md\B(\phi\vee\psi)$. 

        For the converse direction, we assume that $\M,s,t\not\md\B(\phi\vee\psi)$. Then, there is some state $s'\in R(s)$ s.t. $\M,s',t\not\md\phi\vee\psi$, which entails $\M,s,t\not\md\B\phi$ and $\M,s',t\not\md\psi$. Using item (2), we can infer $\M,s,t\not\md\psi$ from $\M,s',t\not\md\psi$. Hence, $\M,s,t\not\md\B\phi\vee\psi$. The proof is completed. \qed
    \end{proof}


\end{document}